\documentclass[11pt,twoside]{amsart}  

\usepackage{amssymb}
\usepackage{latexsym}
\usepackage{epsfig}
\usepackage{psfrag}

\newcommand{\loglike}[1]{\mathop{\rm #1}\nolimits}
\newcommand{\supp}{\loglike{supp}}
\newcommand{\lip}{\loglike{lip}}
\newcommand{\Lip}{\loglike{Lip}}

\newcommand{\eps}{\varepsilon}

\newcommand{\N}{{\mathbb N}}
\newcommand{\R}{{\mathbb R}}

\theoremstyle{plain}
\newtheorem{thm}{Theorem}[section]

\newtheorem{prop}[thm]{Proposition}

\newtheorem{lemma}[thm]{Lemma}

\theoremstyle{definition}
\newtheorem{definition}[thm]{Definition}

\theoremstyle{remark}

\newtheorem{rem}[thm]{Remark}

\numberwithin{equation}{section}

\newcommand{\reff}[1]{(\ref{#1})}

\begin{document}

\title {Lipschitz spaces and $M$-ideals} 

\author[Heiko Berninger and Dirk Werner]
{ Heiko Berninger and Dirk Werner }

\address{Department of Mathematics, Freie Universit\"at Berlin,
Arnimallee~2--6, \qquad {}\linebreak D-14\,195~Berlin, Germany}

\email{berninger@math.fu-berlin.de, werner@math.fu-berlin.de}

\subjclass[2000]{Primary 46B04; secondary 46B20, 46E15}

\keywords{3-ball property, H\"older spaces, Lipschitz spaces, little Lipschitz spaces, $M$-embedded spaces, $M$-ideals}

\begin{abstract}
For a metric space $(K,d)$ the Banach space $\Lip(K)$ consists of all
scalar-valued bounded Lipschitz functions on $K$ with the norm
$\|f\|_{L}=\max(\|f\|_{\infty},L(f))$, where $L(f)$ is the Lipschitz
constant of $f$. The closed subspace $\lip(K)$ of $\Lip(K)$ contains
all elements of $\Lip(K)$ satisfying the $\lip$-condition
$\lim_{0<d(x,y)\to 0}|f(x)-f(y)|/d(x,y)=0$. For
$K=([0,1],|\,{\cdot}\,|^{\alpha})$, $0<\alpha<1$, we prove that
$\lip(K)$ is a proper $M$-ideal in 
a certain subspace of $\Lip(K)$ containing a copy of $\ell^{\infty}$.
\end{abstract}

\maketitle


\section{Introduction}
\label{1}

A scalar-valued function $f$ on a metric space $(K,d)$ satisfying 
$$
L(f)=\sup_{\substack{x,y\in K \\ x\neq y}} \frac{|f(x)-f(y)|}{d(x,y)} < \infty
$$
is called a \emph{Lipschitz function}. The Banach space $\Lip(K)$ of all bounded Lipschitz functions $f$ on $K$ with the norm 
$$
\|f\|_{L}=\max(\|f\|_{\infty},L(f))
$$
is the \emph{Lipschitz space} on $K$. The closed subspace $\lip(K)$ of $\Lip(K)$ which contains all elements of $\Lip(K)$ satisfying the \emph{$\lip$-condition} 
\begin{equation}
\label{ulip}
\lim_{\substack{d(x,y)\to 0 \\ x\neq y}} \frac{|f(x)-f(y)|}{d(x,y)}=0 
\end{equation} 
is called \emph{little Lipschitz space}; its elements are \emph{little Lipschitz functions}. Note that $\Lip(K)$ with the norm $\|f\|_{A}=\|f\|_{\infty}+L(f)$ is a Banach algebra (the \emph{Lipschitz algebra}, see \cite{BCD}), which is not true for the norm $\|\,{\cdot}\,\|_{L}$ (used e.g.\ in \cite{BFT}, \cite{Han92}, \cite{Jenk}, \cite{John70}, and \cite{dL}). Some authors (see \cite{Cies60} and \cite{NW99}) consider the set $\Lip_{0}(K)$ (or $\lip_{0}(K)$) of all (little) Lipschitz functions (possibly unbounded) satisfying $f(x_{0})=0$ for a base point $x_{0}\in K$. Then $L(\cdot)$ is a norm. Equivalently (showing that $\Lip_{0}(K)$ is independent of the choice of $x_{0}$) one can consider $L(\cdot)$ on the factor space $\Lip(K)/N$, where $N$ is the space of the constant functions on $K$ (see \cite{KP} and \cite{Wulb}). 

In the articles mentioned above much work has been done investigating the relationship between $\Lip(K)$ resp.\ $\lip(K)$ and the sequence spaces $\ell^{\infty}$ resp.\ $c_{0}$. In fact the most significant results can be obtained for spaces of \emph{H\"older functions}, i.e.\ Lipschitz functions on $K$ equipped with a H\"older metric $d^{\alpha}$ (where $d^{\alpha}(x,y)=(d(x,y))^{\alpha}$) and $0<\alpha<1$. (We use the notation $K^{\alpha}=(K,d^{\alpha})$ and, unless otherwise stated, we assume $0<\alpha<1$ throughout the paper.) One reason is that in these cases $\lip(K^{\alpha})$ is a nontrivial ``rich'' subspace of $\Lip(K^{\alpha})$. Note that for an open and connected set $K\subset\R^{n}$ with the Eucledian metric $\lip(K)=N$ is one-dimensional whereas for any $K$ the space $\lip(K^{\alpha})$ contains $\Lip(K)$ as a subspace (which is dense if $K$ is compact, see \cite[Proposition~4]{Han92} and \cite[Corollary~1.5]{NW96}). 

Ciesielski \cite{Cies60} seems to have been the first to find a
relationship between Lipschitz spaces and the sequence spaces
$\ell^{\infty}$ and $c_{0}$. For $K=([0,1],|\,{\cdot}\,|)$ and the
H\"older spaces $H_{\alpha}=\Lip_{0}(K^{\alpha})$ and $H_{\alpha}^{0}=
\lip_{0}(K^{\alpha})$ (with base point $0$) he constructs an
isomorphism $T: \ell^{\infty}\to H_{\alpha}$ which carries $c_{0}$
onto $H_{\alpha}^{0}$. The construction uses the well-known Schauder
basis $(\varphi_{n})_{n\in\N}$ for continuous functions on $[0,1]$
(vanishing at 0) which is normalized in $H_{\alpha}$. Via the
uniformly convergent series $\sum_{n=1}^{\infty}a_{n}\varphi_{n}$ one
recovers all elements of $H_{\alpha}$ (resp.\ $H_{\alpha}^{0}$) if the
sequence $(a_{n})_{n\in\N}$ is in $\ell^{\infty}$ (resp.\
$c_{0}$). Adapting this idea in \cite{BFT}, Bonic, Frampton and Tromba
(see also \cite[p.~99]{NW99}) construct an isomorphism
$T:\ell^{\infty}\to \Lip(K^{\alpha})$ with $T(c_{0})=\lip(K^{\alpha})$
for any nontrivial simplex $K\subset\R^{n}$ (with the Eucledian
metric). Although their result is also true for any finite dimensional
(and infinite) compact set $K$, it is still an open question 
whether this can be carried over to arbitrary compact sets $K$.

However, for any metric space $K$ with $\Lip(K)\neq\lip(K)$, thus with $\inf_{d(x,y)>0}d(x,y)=0$, Johnson (cf.\ \cite{John72} and \cite{John74}) gives an isomorphism $T:\ell^{\infty}\to T(\ell^{\infty})\subset \Lip(K^{\alpha})$ for which $T(c_{0})$ is complemented in $\lip(K^{\alpha})$. Consequently $\lip(K^{\alpha})$ is not complemented in $\Lip(K^{\alpha})$ for these metric spaces. In general $\Lip(K)$ contains a copy of $\ell^{\infty}$ unless it is finite dimensional. 

The isometric representation of the spaces $\Lip_{0}(K^{\alpha})$ and
$\lip_{0}(K^{\alpha})$ with $\alpha>0$ and compact $K$ is the subject
of Wulbert's article \cite{Wulb}. His stunning result is (provided a
gap in the proof of his Proposition 3.1 can be filled) that a point
separating little Lipschitz space $\lip_{0}(K^{\alpha})$ can only be
isometrically isomorphic to $c_{0}$ if $\alpha=1$ and $K$ is isometric
to a nowhere dense subset of the real line. In turn such an isometric
isomorphism exists provided $K$ is in addition a Lebesgue null set
(see also \cite[pp.\ 73--76]{B}). This isometry is easily obtained by
encoding the slopes of certain chords of a function on $K$
 in a sequence.

Due to a duality result first observed by de Leeuw in \cite{dL} it is
 possible to trace the problem of an (isometric) isomorphism between
 $\Lip(K)$ and $\ell^{\infty}$ to the same problem for $\lip(K)$ and
 $c_{0}$. For the circle $K\subset\R^{2}$ with the metric given by
 unit arc length, he proves that $I:\lip(K^{\alpha})''\to
 \Lip(K^{\alpha})$ defined by $I(F)(x)=F(\delta_{x})$ (where
 $\delta_{x}$ is a point evaluation functional) is an isometric
 isomorphism. This isometry is natural in the sense that $I\circ\pi$
 is the identity on $\lip(K^{\alpha})$ if $\pi$ is the natural
 embedding of $\lip(K^{\alpha})$ into
 $\lip(K^{\alpha})''$. One key idea for the proof is to identify $\lip(K^{\alpha})$ as a subspace of a space of continuous functions vanishing at infinity on a certain locally compact space and then to use the Riesz representation theorem. It was Jenkins in \cite{Jenk} who refined de Leeuw's technique to get the result for any compact (and certain locally compact) H\"older metric spaces $K^{\alpha}$ which (for spaces of complex-valued functions) satisfy an additional geometric condition. Jenkins' idea is to approximate $\Lip(K^{\alpha})$-functions for $\alpha<\beta\leq 1$ uniformly by $\Lip(K^{\beta})$-functions with only slightly increased norm (note that $\Lip(K^{\beta})\subset \lip(K^{\alpha})$) by applying an extension theorem (see \cite{McS}) for Lipschitz functions. 

In \cite[p.~159]{John70} Johnson extracts this essential property (S)
of H\"older functions from Jenkins' key lemma and obtains the result
by considering the pre-adjoint of $I$. In addition he gets rid of the
geometric condition Jenkins imposed on $K$ for the complex case
showing that it was only technical in nature. Apparently only aware of
de Leeuw's article Bade, Curtis and Dales rediscovered in \cite{BCD}
the general result (for Lipschitz algebras) pointing out the duality
$\langle f,\mu\rangle$, $f\in \Lip(K)$, $\mu\in\lip(K)'$, and giving
the inverse of $I$ if condition (S) is fulfilled (in fact they use an
equivalent condition which seems to be weaker than (S)). Hanin pursues
yet another approach in \cite{Han92} and \cite{Han94} by identifying
(for real-valued functions) $\lip(K)'$ with the closure of the space
of Borel measures on $K$ (which is assumed to be compact or to have
compact closed balls, see \cite{Han94}) equipped with an appropriate
norm, the Kantorovich-Rubinstein norm (see
\cite[pp. 225--237]{Kant}). Moreover, but again in the real case, he
proves that (S) characterizes all compact metric spaces for which $I$
is an isometric isomorphism. In \cite[p.~287]{NW96} and
\cite[p.~2644]{NW97} Weaver shows that this equivalence also applies in
the case of complex-valued functions and gives (with a slight error in
\cite[Theorem 1.(d)]{NW97}) different equivalent versions of (S),
which can be interpreted as a kind of uniform point separation
property of $\lip(K)$. Adapting the proof in \cite{BCD} he extends in
\cite{NW96a} the characterization result to a certain class of locally
compact metric spaces. For a comprehensive overview of the facts known
so far about Lipschitz spaces we also refer 
to \cite{NW99}.
            
The isomorphism results $\Lip(K)\simeq\ell^{\infty}$, $\lip(K)\simeq
c_{0}$ as well as the duality $\lip(K)''\cong \Lip(K)$ suggest a close
relationship between the Lipschitz spaces $\Lip(K)$ resp.\ $\lip(K)$
and the sequence spaces $\ell^{\infty}$ resp.\ $c_{0}$ at least for
(infinite) compact H\"older metric spaces $K\subset \R^{n}$. The
duality result still holds for metric spaces $K$ equipped with
generalized H\"older metrics (see \cite[p.~348]{Han94}), however it
fails in the simple finite-dimensional cases $K=[0,1]^n$ with the
Eucledian metric since then $\lip(K)$ is trivial. And although the
isomorphism result $\Lip(K)\simeq\ell^{\infty}$ is still true for
(infinite) compact sets $K\subset\R$ (see \cite[1.1]{MuSem} and
\cite[p.~480]{John75}), it is false for $K=[0,1]^2$ (which can be
shown just as in the case of continuously differentiable functions on
$[0,1]^2$, see \cite[p.~59]{MW}). This indicates that H\"older spaces
$\Lip(K^{\alpha}$) resp.\ $\lip(K^{\alpha})$ form a quite special
subclass of the class of all Lipschitz spaces $\Lip(K)$ resp.\
$\lip(K)$. 

In the following we will have these spaces and their connection with
the sequence spaces $\ell^{\infty}$ and $c_{0}$ in mind when studying
$\Lip(K)$ and $\lip(K)$ with regard to their $M$-structure. 
In Section~\ref{2} of this paper we shall prove our main theorem which
states that
$H_{\alpha}^{0}$ is a
proper $M$-ideal in a certain nonseparable subspace of
$H_{\alpha}$ (for the definition of an
$M$-ideal see below). 
On the way to this theorem we
first state a lemma which applies in a general setting and which
contains our main idea of approximating big Lipschitz functions by
little Lipschitz functions in order to check if the $3$-ball property
(cf.\ Proposition~\ref{3KP})
is fulfilled (Lemma~\ref{3B}). Then, after analysing what it means
that a function is ``little Lipschitz'' at a point 
(Definition~\ref{unif lip}), 
we give the definition of ``pointwise big Lipschitz''
functions and collect them in the space $H_{\alpha}^{p}$ 
(Definition~\ref{Hap}). We prove (by ``inserting constants'') that
$H_{\alpha}^{0}$ satisfies the $3$-ball property in $H_{\alpha}^{p}$
(Proposition~\ref{Ha0 3B Hap}) and by a polygonal construction
characterize the closure of $H_{\alpha}^{p}$ in $H_{\alpha}$ as the
space $H_{\alpha}^{\omega}$ of ``weakly big Lipschitz'' functions
(Definition~\ref{Haw} and Proposition~\ref{Hap dense Haw}), in which
$H_{\alpha}^{0}$ is an $M$-ideal. Finally by an application of the
mean value theorem we show that $H_{\alpha}^{0}$ is a proper $M$-ideal
in $H_{\alpha}^{\omega}$ (Proposition~\ref{Ha0 pMid Haw}). Using a
more general construction of Johnson in \cite{John74} it becomes clear
that $H_{\alpha}^{0}$ is not even complemented in
$H_{\alpha}^{\omega}$ (see Remark~\ref{Mid-rem}).
Our  results will be recast in terms of the
Ciesielski basis of $H^0_\alpha$ in Section~\ref{3}. It remains open,
however, whether  $H^0_\alpha$ is an $M$-ideal in $H_\alpha$; in
Section~\ref{4} we are going to study a certain Cantor-like $H_\alpha$-function
that is very remote from $H^0_\alpha$ and that might lead to a
negative answer to this question.

\bigskip\noindent
\textbf{Acknowledgement.}
The authors thank P.~Wojtaszczyk for fruitful discussions on
the subject of this note.


\section{$M$-ideals of Lipschitz functions}
\label{2}

We start by recalling the definition of an $M$-ideal, introduced in
  \cite{AE}.  A detailed study of this notion can be found in \cite{HWW}.

\begin{definition}
\label{M-ideal}
Let $X$ be a real or complex Banach space.
\begin{enumerate}
\item
A linear projection $P$ on $X$ is called an \emph{$M$-projection} if
\begin{equation*}
\|x\|=\max(\|P(x)\|,\|x-P(x)\|)\quad\forall x\in X
\end{equation*}
and an \emph{$L$-projection} if
$$
\|x\|=\|P(x)\|+\|x-P(x)\|\quad\forall x\in X.
$$
\item
A closed subspace $U\subset X$ is called an \emph{$M$-summand} resp.\ 
\emph{$L$-summand} if it is the range of an $M$-projection resp.\ 
$L$-projection.
\item
A closed subspace $U\subset X$ is called an \emph{$M$-ideal} if its
annihilator $U^{\perp}$ is an $L$-summand in $X'$. It is called a
\emph{proper} $M$-ideal if it is not an $M$-summand. If $U$ is an
$M$-ideal in its bidual, then $U$ is 
called an \emph{$M$-embedded} space.
\end{enumerate}
\end{definition}

It is well known and easy to see that $c_{0}$ is an $M$-embedded
space. Moreover it is a kind of prototype of a proper $M$-ideal since
any proper $M$-ideal $U$ contains a copy of $c_{0}$ (which is
complemented if $U$ is $M$-embedded), see \cite[II.4.7 and
III.4.7]{HWW}. Now, in light of the results given in the introduction, it
is natural to ask whether (under reasonable conditions on $K$) $\lip(K)$
is an $M$-ideal in $\Lip(K)$. So far this is only clear in the very
restricted isometric cases given by Wulbert in \cite[Lemma~3.4]{Wulb}
where $K$ is a nowhere dense (infinite) compact subset of $\R$ with
Lebesgue measure $0$. On the other hand $\lip(K)$ is not an $M$-ideal
in the Lipschitz algebra $(\Lip(K),\|\,{\cdot}\,\|_{A})$ (provided
$\lip(K)\neq \Lip(K)$) since in a unital commutative Banach algebra an
$M$-ideal is always a closed ideal (see \cite[V.4.1]{HWW}). Therefore
we will from now on consider the Lipschitz spaces
$(\Lip(K),\|\,{\cdot}\,\|_{L})$ and $(\Lip_{0}(K),L(\cdot))$. In
\cite[Lemma~3.4]{BCD} (although in that paper Lipschitz algebras are
considered) one gets a hint of what a natural projection
$P:\Lip(K)'\to\lip(K)^{\perp}$ could look like. Nevertheless if one
tries to check the norm equation given in the definition one always
ends up examining the relationship between the elements of the closed
unit sphere $B_{U}$ of $U$ and the elements in $B_{X}$. This fact is
reflected by the following well-known equivalence (see \cite[I.2.2]{HWW}).

\begin{prop}
\label{3KP}
A closed subspace $U$ of a Banach space $X$ is an $M$-ideal in $X$ if and only if for all $y_{1},y_{2},y_{3}\in B_{U}$, all $x\in B_{X}$ and all $\eps>0$ there is some $y\in U$ satisfying
$$
\|x+y_{i}-y\|\leq 1+\eps\quad (i=1,2,3).
$$
\end{prop}

This condition is called the (restricted) $3$-ball property of $M$-ideals. We will also say that an element $x\in B_{X}$ satisfies the $3$-ball property if the condition is fulfilled for this $x$. 

We are now going to  prove  that
$H_{\alpha}^{0}=\lip_{0}([0,1]^{\alpha})$ (with base point $0$) is a
proper $M$-ideal in a certain nonseparable subspace of
$H_{\alpha}=\Lip_{0}([0,1]^{\alpha})$. 
There is a standard technique how to verify the $3$-ball property in
concrete function or sequence spaces (mostly in the case of
$M$-embedded spaces) which can be found in \cite[p.~102]{HWW} and in
\cite{DW92}. In many cases in which a Banach space $U$ might be an
$M$-ideal in a second space $X$ (possibly its bidual) one can observe
an $o(\cdot)$-$O(\cdot)$-relation between the elements of $U$ and the
elements of $X$, where $o$ and $O$ denote the Landau symbols. Examples
of this are of course $c_{0}$ and $\ell^{\infty}$ or $\lip(K)$ and
$\Lip(K)$. Suppose the elements of $X$ are functions on a set $K$ for
which a value $|\,{\cdot}\,|$ is bounded (the $O$-condition). In order
to verify the $3$-ball property for given $y_{1},y_{2},y_{3}\in
B_{U}$, $x\in B_{X}$ and $\eps>0$, choose a subset $M\subset K$ such
that $|y_{i}|\leq\eps$, $i=1,2,3$, on $K\backslash M$ (the
$o$-condition). Then try to define $y\in U$ by $y\approx x$ up to
$\eps>0$ on $M$ and elsewhere not ``too far'' from $x$ such that
$|x-y|\leq 1+\eps$ holds on $K\backslash M$. Then on $M$ one gets
$|x+y_{i}-y|\leq |y_{i}|+|x-y|\leq 1+\eps$ while on $K\backslash M$
the $o$-condition leads to $|x+y_{i}-y|\leq |x-y|+|y_{i}|\leq
1+2\eps$. Observe that although the $o$-condition for $y$ leads to
$|y|\leq\eps$ on a certain set $M'\subset K$ providing $|x-y|\leq
1+\eps$ on $K\backslash M'$ this might not help much since $M'$
depends on $M$. On the other hand it is obviously enough to check the
$3$-ball property for elements in 
dense subsets of $B_{U}$ or $B_{X}$ respectively.

The above method can be perfectly applied to the standard example
$c_{0}$ in $\ell^{\infty}$. In the case of $\lip(K)$ in $\Lip(K)$ it
is less obvious how to extract a ``set'' $M\subset K$ since the
$o$-condition in $\lip(K)$ involves all elements of $K$, however
locally. Thus we seek for a local condition to replace the set
$K\backslash M$. The following abbreviation will be helpful 
in the sequel.

\begin{definition}
\label{slope}
Let $f$ be a Lipschitz function on a metric space $K$ and $x,y\in K$, $x\neq y$. Then 
$$
L_{xy}(f)=\frac{|f(x)-f(y)|}{d(x,y)}
$$
defines the \emph{slope} of $f$ between $x$ and $y$.
\end{definition}

As a motivation we give a version of the separation property (S) mentioned in the introduction which for compact metric spaces $K$ is sufficient and necessary for the duality result (cf.\ \cite[Theorem~1]{NW97} and \cite[Corollary 3.3.5]{NW99}):

\begin{itemize}
\item[(S)]  
\emph{For any constant $c>1$, any finite set $A\subset K$ and any function $h\in B_{\Lip(K)}$ there is a function $g\in\lip(K)$ satisfying $g_{|A}=h_{|A}$ and $\|g\|_{L}\leq c$.}
\end{itemize}

Analogous versions exist for $\Lip_{0}$-spaces. In particular the separation property, which is satisfied for all compact H\"older metric spaces $K^{\alpha}$, allows uniform approximation of big Lipschitz functions $h$ by little Lipschitz functions $g$. Using that $L_{xy}(h-g)$ becomes small for small $\|h-g\|_{\infty}$ provided $d(x,y)$ is bounded below we obtain the local condition sought for which supplies a criterion for the $3$-ball property of little Lipschitz spaces.

\begin{lemma}
\label{3B}
Let $K$ be a metric space and $W$ a closed subspace of $\Lip(K)$ with 
$\lip(K)\subset W$ satisfying the following condition \emph{(3B):}

For any $h\in B_{W}$, $\eps'>0$ and  $\delta'>0$ there exists
some $g\in \lip(K)$ with
\begin{equation}
\label{3B1}
\|h-g\|_{\infty}\leq\eps'\delta'
\end{equation}
and
\begin{equation}
\label{3B2}
0<d(x,y)\leq\delta'\Longrightarrow L_{xy}(h-g)\leq 1+\eps'\quad\forall x,y\in K.
\end{equation}
Then $\lip(K)$ is an $M$-ideal in $W$.

The analogous statement holds true in $\Lip_{0}(K)$.
\end{lemma}

\begin{proof}
We show that $\lip(K)$ satisfies the $3$-ball property in $W$. So let $h\in B_{W}$, $f_{1},f_{2},f_{3}\in B_{\lip(K)}$ and $\eps'>0$. Using the $\lip$-condition choose $\delta'>0$ such that $0<d(x,y)\leq\delta'$ forces
$$
L_{xy}(f_{i})=\frac{|f_{i}(x)-f_{i}(y)|}{d(x,y)}\leq\frac{\eps'}{2}\quad\forall i\in\{1,2,3\}.
$$
Now condition (3B) and~\reff{3B1} yield some $g\in \lip(K)$ with $\|h-g\|_{\infty}\leq\delta'\eps'/4$. Thus for $d(x,y)\geq\delta'$ we obtain
$$
L_{xy}(h-g)\leq \frac{|(h-g)(x)|+|(h-g)(y)|}{d(x,y)}\leq\frac{\eps'}{2}.
$$
Altogether we conclude
$$
L_{xy}(h-g+f_{i})\leq L_{xy}(h-g)+L_{xy}(f_{i})\leq \frac{\eps'}{2}+1\leq 1+\eps'\quad\forall i\in\{1,2,3\}
$$
if $d(x,y)\geq \delta'$ while 
$$
L_{xy}(h-g+f_{i})\leq L_{xy}(h-g)+L_{xy}(f_{i})\leq 1+\frac{\eps'}{4}+\frac{\eps'}{2}\leq 1+\eps'\quad\forall i\in\{1,2,3\}
$$
holds true for $0<d(x,y)\leq \delta'$ due to~\reff{3B2} and the choice
of $\delta'$.

The inequality
$$
\|h-g+f_{i}\|_{\infty}\leq\|h-g\|_{\infty}+\|f_{i}\|_{\infty}\leq \eps'+1
$$
is obtained for any $\delta'\leq 1$ according to~\reff{3B1}. 

The analogous statement is obvious in $\Lip_{0}(K)$.
\end{proof}

For short, condition (3B) requires uniform approximation of a big Lip\-schitz function $h$ by a little Lipschitz function $g$ such that $g$ ``adapts the slope of $h$ locally''.

We now turn to the application of Lemma~\ref{3B} to  the H\"older
spaces $H_{\alpha}$ and $H_{\alpha}^{0}$ on the unit interval. We
still assume $0<\alpha<1$ unless otherwise stated and it should not
lead to confusion if both the Lipschitz constant $L(f)$ and
$L_{xy}(f)$ are now used for the underlying metric
$d(x,y)=|x-y|^{\alpha}$ 
so that $L(f)$ is actually the H\"older constant of $f$. 

One succeeds to check the $3$-ball property applying Lemma~\ref{3B}
to quite simple functions in $H_{\alpha}$ like $h:x\mapsto
x^{\alpha}$ (or functions built from ``$x^{\alpha}$-arcs'') when one
approximates them by lines around the ``critical points''. More
generally consider a partition $0=x_{0}<x_{1}<\dots<x_{n}=1$ of the
unit interval and a continuous function $p$ on $[0,1]$ vanishing at
$0$ and affine on $[x_{k-1},x_{k}]$ for $k=1,\dots,n$. We say that $p$
is a \emph{polygon} or a \emph{polygonal function} and call the points
$(x_{k},p(x_{k}))$ the \emph{nodes} of $p$. It is obvious that $p\in
H_{1}\subset H_{\alpha}^{0}$. As indicated in the introduction (see
\cite{Cies60} or \cite[1.5]{MuSem}) the set of all polygonal functions
is even dense in $H_{\alpha}^{0}$. As an important tool we state a
result of Krein and Petuin in 
\cite[Lemma~5.1]{KP}. 

\begin{lemma}
\label{KP}
Let $h\in H_{\alpha}$ and $p$ be a polygon interpolating $h$ in its nodes $(x_{k},h(x_{k}))$, $k=1,\dots,n$. Then $L_{xy}(p)\leq L_{x_{k-1},x_{k}}(h)$ if $x,y\in [x_{k-1},x_{k}]$ for some $k\in\{1,\dots,n\}$ and moreover $L(p)\leq L(h)$.
\end{lemma}

Observe that this lemma implies the separation property of
$H_{\alpha}^{0}$ with $c=1$. By polygonal interpolation one can prove
using this lemma that any monotonic $h\in B_{H_{\alpha}}$ satisfies
the $3$-ball property. However by considering strongly oscillating
functions (cf.\ Definition~\ref{almond function}) one can show that
there is no hope that polygonal interpolation (as described in the
above lemma) will succeed in proving 
the 3-ball property with our criterion from Lemma~\ref{3B} for any
$h\in B_{H_{\alpha}}$ just by choosing the partition fine enough. This
is remarkable for two reasons: First of all it demonstrates that a
straightforward application of the separation property fails to
satisfy this quite natural criterion and secondly if $H_{\alpha}^{0}$
is an $M$-ideal in $H_{\alpha}$ some polygonal approximation
\emph{has} to work since the polygons are dense in 
$H_{\alpha}^{0}$. 

Now we restrict our considerations to big Lipschitz functions with
only finitely many ``critical points'', a notion that we formalize
first. Quite naturally one could think of $\limsup_{y\to x}
L_{xy}(h)>0$ as an appropriate condition for a critical point $x$ of a
Lipschitz function $h$. But this turns out to be no good since the
converse $\lim_{y\to x} L_{xy}(h)=0$ is just a pointwise
$\lip$-condition. Except for Krein and Petuin in \cite{KP} all authors
listed in the references use the uniform $\lip$-condition \reff{ulip} and 
for the essential isomorphism and duality results the uniformity is crucial. 

\begin{prop}
\label{pw lip}
Let $0<\alpha<1$.
There is a function $h$ on $[0,1]$ with $h\notin H_{\alpha}$
satisfying the pointwise $\lip$-condition with respect to the metric
$d^\alpha$ at every point. 
\end{prop}

\begin{proof}
Let $x_{k}=2^{-k}$ for $k\in\N$. We shall define $h$ on $[0,1]$ with
$0\leq h(x)\leq x$ for all $x$ such that $\lim_{y\to 0}
L_{0y}(h)=0$. Let $(\delta_{k})_{k\in\N}$ be a positive sequence
satisfying $x_{k+1}+\delta_{k+1}\leq x_{k}-\delta_{k}$ and $x_{k}\leq
k(\delta_{k})^{\alpha}$. 
Now define
$$
h(x_k):=k(\delta_{k})^{\alpha}, \  h(x_k\pm\delta_k):=0
\quad\forall k\in\N,  
$$
choose $h$ linear on $[x_{k}-\delta_{k},x_{k}]$ and
$[x_{k},x_{k}+\delta_{k}]$ 
for any $k\in\N$ and $0$ elsewhere on $[0,1]$.

By construction $h$ satisfies the pointwise $\lip$-condition at~$0$, but
since $h$ is polygonal on $[x,1]$ for any positive $x<1$ it is a
Lipschitz function with exponent~$1$ on $[x,1]$; thus it even
satisfies the uniform $\lip$-condition for any $0<\alpha'<1$ and any
point in $[x,1]$. However $L_{x_k,x_{k}+\delta_k}(h)=k$ 
for all $k\in\N$ so that $h\notin H_{\alpha}$.
\end{proof} 

By a modification of the above argument one can even exhibit
a continuous function $h\notin \bigcup_{\beta\in (0,1)} H_\beta$ that
satisfies the pointwise $\lip$-condition for 
$d^\alpha$ at every point. 

There is a straightforward generalization of Proposition~\ref{pw lip}
to all metric spaces $K^{\alpha}$ possessing infinitely many cluster
points. Now it is sensible enough to call the behaviour of the
function $h$ in the above proof ``critical'' at the point~$0$. Here we
decide on a quite restrictive definition of a noncritical point and we
formulate it for the general case.

\begin{definition}
\label{unif lip}
Let $h$ be a function on a metric space $K$. We say that $h$
\emph{satisfies the $\lip$-condition at the point $x\in K$} and call
$x$ a \emph{noncritical} point of $h$ if there is a neighbourhood
$U_{x}$ of $x$ in which $h$ satisfies the uniform
$\lip$-condition~\reff{ulip}. Otherwise $x$ is called a
\emph{critical} point of~$h$.
\end{definition}

Now in contrast to the counterexample in Proposition~\ref{pw lip} we
can localize 
the little Lipschitz property of a function.

\begin{lemma}
\label{loc lip=lip}
Let $g$ be a function on a compact metric space $K$ with no critical point. Then $g\in\lip(K)$.
\end{lemma}

\begin{proof}
First we note that $g$ satisfies the uniform $\lip$-condition on $K$
since otherwise there exist $\eps>0$ and sequences $(x_{k})$ and
$(y_{k})$ in $K$ both converging to a limit $x\in K$ with
$L_{x_k,y_k}(g)\geq\eps$ and $x$ is a critical point. We now show that
$g$ is a Lipschitz function. Assume not; then again there are
sequences $(x_{k})$ and $(y_{k})$ in $K$ converging to limits $x$ and
$y$ respectively with $L_{x_k,y_k}(g)\to\infty$. Since $g$ satisfies
the (uniform) $\lip$-condition we have $x\neq y$. But then since
$L_{x_k,y_k}(g)$ is unbounded $g$ is also unbounded contradicting the
continuity of $g$ implied by the 
$\lip$-condition.
\end{proof}

Now we return to the one-dimensional situation.

\begin{definition}
\label{Hap}
By $H_{\alpha}^{p}$ we denote the subspace of all functions in
$H_{\alpha}$ 
which have finitely many critical points.
\end{definition}

Applying Lemma~\ref{3B} we now prove that $H_{\alpha}^{0}$ satisfies
the 3-ball property in $H_{\alpha}^{p}$ by approximating an
$H_{\alpha}^{p}$-function $h$ in the following way: Keep the function
untouched outside small neighbourhoods $(x_{k}-\delta,x_{k}+\delta)$
around the critical points $x_{k}$ of $h$ and constant inside them. By
adding appropriate constants on the different intervals 
we obtain a function $g\in H_{\alpha}^{0}$ which does the job.

\begin{prop}
\label{Ha0 3B Hap}
$H_{\alpha}^{0}$ satisfies the 3-ball property in $H_{\alpha}^{p}$.
\end{prop}

\begin{proof}
Let $h\in B_{H_{\alpha}^{p}}$ with its finitely many critical points
$x_{1}<\dots<x_{n}$ and for simplicity include without loss of
generality $x_{1}=0$ and $x_{n}=1$. Let $\eps'>0$ and $\delta'>0$
where without loss of generality $(\delta')^{1/\alpha}\leq
\frac{1}{2}\min_{1\leq k\leq n-1}|x_{k}-x_{k+1}|$, and set
$\eps:=\eps'\delta'$. Now choose a positive
$\delta\leq\min(\frac{1}{2}(\delta')^{1/\alpha},\frac{1}{2}(\frac{\eps}{n})^{1/\alpha})$
such that $|x-y|\leq 2\delta$ 
implies $|h(x)-h(y)|\leq\frac{\eps}{n}$. We set $h(-\delta):=0$,
$$
\varDelta h(x_{k}):=h(x_{k}-\delta)-h(x_{k}+\delta)\quad\forall k\in \{1,\dots,n\}
$$
and define the function $g$ on $[0,1]$ by
$$
g(x):=\left
\{\begin{array}{ll}
h(x_{k}-\delta)+\sum_{i=1}^{k-1}\varDelta h(x_{i}) & \mbox{ if } x_{k}-\delta\leq x\leq x_{k}+\delta\\
h(x)+\sum_{i=1}^{k}\varDelta h(x_{i}) & \mbox{ if } x_{k}+\delta\leq x\leq x_{k+1}-\delta.
\end{array}\right.
$$
By construction $g$ is continuous and $g(0)=0$. We show $\|h-g\|_{\infty}\leq\eps$. If $x_{k}-\delta\leq x\leq x_{k}+\delta$, $k\in\{1,\dots,n\}$, our choice of $\delta$ provides
$$
|h(x)-g(x)|\leq |h(x)-h(x_{k}-\delta)|+ \sum_{i=1}^{k-1}|\varDelta h(x_{i})|\leq\frac{\eps}{n}+(k-1)\frac{\eps}{n}\leq\eps
$$
while for $x_{k}+\delta\leq x\leq x_{k+1}-\delta$, $k=1,\dots,n-1$, we get
$$
|h(x)-g(x)|\leq\sum_{i=1}^{k}|\varDelta h(x_{i})|\leq(n-1)\frac{\eps}{n}\leq\eps.
$$
Now we prove that $g$ satisfies the $\lip$-condition at every point
$x\in [0,1]$, whence $g\in H_{\alpha}^{0}$ by Lemma~\ref{loc
  lip=lip}. If $x$ is not of the form $x_{k}\pm\delta$ for a
$k\in\{1,\dots,n\}$, then in a neighbourhood of~$x$ either $g$ is
constant or it differs from $h$ by a constant. In both cases $x$ is a
noncritical point of $g$ (observe that in the latter case $x$ is a
noncritical point of $h$). We still have to estimate $L_{xy}(g)$ for
$x,y$ in a neighbourhood of a point $x_{k}-\delta$ or
$x_{k}+\delta$. If $x_{k-1}+\delta< x<x_{k}-\delta<y<x_{k}+\delta$ 
we have
\begin{align*}
L_{xy}(g) & =\frac{\left|\left(h(x)+\sum_{i=1}^{k-1}\varDelta h(x_{i})\right)-\left(h(x_{k}-\delta)+\sum_{i=1}^{k-1}\varDelta h(x_{i})\right)\right|}{|x-y|^{\alpha}}\\
& \leq\frac{|h(x)-h(x_{k}-\delta)|}{|x-(x_{k}-\delta)|^{\alpha}}=L_{x,x_{k}-\delta}(h)
\end{align*}
and quite similarly if $x_{k}-\delta< x<x_{k}+\delta <y<x_{k+1}-\delta$  
\begin{align*}
L_{xy}(g)\leq L_{x_{k}+\delta,y}(h).
\end{align*}  
In the remaining cases in which both $x$ and $y$ are (locally) above
or below $x_{k}-\delta$ or $x_{k}+\delta$ respectively we clearly have
$L_{xy}(g)\leq L_{xy}(h)$. Altogether these estimates show that since
$h$ satisfies the $\lip$-condition at the points $x_{k}-\delta$,
$k=2,\dots,n$, and $x_{k}+\delta$, $k=1,\dots,n-1$, 
so does $g$.

Finally we prove that $0<|x-y|^{\alpha}\leq\delta'$ implies
$L_{xy}(h-g)\leq 1$. Due to
$\delta\leq\frac{1}{2}(\delta')^{1/\alpha}\leq \frac{1}{4}\min_{1\leq
  k\leq n-1}|x_{k}-x_{k+1}|$ one of the following five cases applies. 
We skip the calculations since they are carried out as above. 

\vspace{0.3cm}

1) $x_{k-1}+\delta\leq x<y\leq x_{k}-\delta$
\begin{equation*}
\Longrightarrow L_{xy}(h-g)=L_{xy}\left(\sum_{i=1}^{k-1}\varDelta h(x_{i})\right)=0.
\end{equation*}

2) $x_{k-1}+\delta< x\leq x_{k}-\delta<y\leq x_{k}+\delta$
\begin{align*}
\Longrightarrow L_{xy}(h-g)\leq L_{x_{k}-\delta,y}(h)\leq 1.
\end{align*}

3) $x_{k}-\delta\leq x<y\leq x_{k}+\delta$
\begin{equation*}
\Longrightarrow L_{xy}(h-g)=L_{xy}(h)\leq 1, \mbox{ since $g=\mbox{const.}$ on $[x_{k}-\delta,x_{k}+\delta]$.}
\end{equation*}

4) $x_{k}-\delta\leq x< x_{k}+\delta\leq y$
\begin{align*}
\Longrightarrow L_{xy}(h-g)\leq L_{x,x_{k}+\delta}(h)\leq 1.
\end{align*}

5) $x_{k-1}+\delta< x\leq x_{k}-\delta<x_{k}+\delta\leq y< x_{k+1}-\delta$
\begin{align*}
\Longrightarrow L_{xy}(h-g)\leq L_{x_{k}-\delta,x_{k}+\delta}(h)\leq 1.
\end{align*}

The assertion now follows from Lemma~\ref{3B}.
\end{proof}
 
It is easy to see, e.g.\ by considering absolutely convergent series $\sum_{k=1}^{\infty}h_{k}$ of functions $h_{k}\in H_{\alpha}^{p}$ with increasing numbers of critical points, that $H_{\alpha}^{p}$ is not closed in $H_{\alpha}$.

\begin{definition}
\label{Haw}
We define $H_{\alpha}^{\omega}$ as the subspace of all functions $h\in H_{\alpha}$ having the following property:

For any $\eps>0$ there are finitely many points $x_{1},\dots,x_{n}$ in $[0,1]$ such that for any $\tilde{x}\in [0,1]\backslash\{x_{k}\}_{k=1}^{n}$ there is a neighbourhood $U(\tilde{x})$ with
$$
\sup_{\substack{x,y\in U(\tilde{x}) \\ x\neq y}} \frac{|h(x)-h(y)|}{|x-y|^{\alpha}}\leq\eps.
$$
\end{definition}

\begin{prop}
\label{Hap dense Haw}
$H_{\alpha}^{\omega}$ is the closure of $H_{\alpha}^{p}$ in $H_{\alpha}$.
\end{prop}

\begin{proof}
We first prove that $H_{\alpha}^{\omega}$ is closed in $H_{\alpha}$. Let $(h_{i})_{i\in\N}$ be a sequence in $H_{\alpha}^{\omega}$ with a limit $h\in H_{\alpha}$ and $\eps>0$. Choose $m\in \N$ such that $L(h_{m}-h)\leq\frac{\eps}{2}$. There are finitely many points $x_{1},\dots,x_{n}$ such that $\sup_{x,y\in U(\tilde{x})}L_{xy}(h_{m})\leq\frac{\eps}{2}$ holds for all $\tilde{x}\in [0,1]\backslash\{x_{k}\}_{k=1}^{n}$ and certain neighbourhoods $U(\tilde{x})$. We obtain
$$
L_{xy}(h)\leq L_{xy}(h_{m})+L_{xy}(h_{m}-h)\leq\frac{\eps}{2}+\frac{\eps}{2}=\eps
$$
for $x,y\in U(\tilde{x})$ and thus $h\in H_{\alpha}^{\omega}$.

Now let $h\in B_{H_{\alpha}^{\omega}}$ and $\eps>0$. Find finitely many points $x_{1}<\dots<x_{n}$ in $[0,1]$ (without loss of generality $x_{1}=0$ and $x_{n}=1$) such that for all $\tilde{x}\in [0,1]\backslash\{x_{k}\}_{k=1}^{n}$ we have $\sup_{x,y\in U(\tilde{x})}L_{xy}(h)\leq\frac{\eps}{4}$ in certain neighbourhoods $U(\tilde{x})$, which we may assume to be open intervals. In the following we define a function $f\in H_{\alpha}^{p}$ approximating $h$ in $H_{\alpha}$ by a polygonal construction around the critical points $x_{1}<\dots<x_{n}$. Let $(\delta_{m})_{m=0}^{\infty}$ be a strictly decreasing positive null sequence with $\delta_{0}<\frac{1}{2}\min_{k=1,\dots,n-1}|x_{k}-x_{k+1}|$. 

From the union of the above sets $U(\tilde{x})$ covering the compact
set $M_{0}:=[0,1]\backslash
\bigcup_{k=1}^{n}(x_{k}-\delta_{0},x_{k}+\delta_{0})$ we obtain a
finite covering $\bigcup_{j=1}^{N_{0}}U(\tilde{x}_{j})$. Define
$\tilde{\ell}_{0}$ as half of the minimal length of all nonempty
intersections $U(\tilde{x}_{i})\cap U(\tilde{x}_{j})$, $1\leq i,j\leq
N_{0}$. Now choose 
$n_{0}^{k}\in \N$ such that 
$$
\frac{|(x_{k+1}-\delta_{0})-(x_{k}+\delta_{0})|}{n_{0}^{k}}=:\ell_{0}^{k}\leq\tilde{\ell}_{0}\quad\forall k=1,\dots,n-1
$$
and define $f$ on $M_{0}$ by
$$
f(x):=\left\{
\begin{array}{ll}
h(x_{k}+\delta_{0}+j\ell_{0}^{k})\  & \mbox{if }\, x=x_{k}+\delta_{0}+j\ell_{0}^{k}\\
& 0\leq j\leq n_{0}^{k},\; k=1,\dots,n-1\\
\mbox{affine on} & [x_{k}+\delta_{0}+j\ell_{0}^{k}, x_{k}+\delta_{0}+(j+1)\ell_{0}^{k}]\\
& 0\leq j\leq n_{0}^{k}-1,\; k=1,\dots,n-1.
\end{array}
\right.
$$
For $m\in\N$ we define $f$ on the compact set
$$
M_{m}:=\left(\bigcup_{k=1}^{n}[x_{k}-\delta_{m-1},x_{k}+\delta_{m-1}]\backslash (x_{k}-\delta_{m},x_{k}+\delta_{m})\right)\,\cap\enspace [0,1]
$$
in the following manner. Again consider a finite covering $\bigcup_{j=N_{m-1}+1}^{N_{m}}U(\tilde{x}_{j})$ of $M_{m}$ obtained from the union of the sets $U(\tilde{x})$. Again let $\tilde{\ell}_{m}$ be half of the minimal length of all nonempty intersections $U(\tilde{x}_{i})\cap U(\tilde{x}_{j})$, $N_{m-1}+1\leq i,j\leq N_{m}$. With a certain $n_{m}\in \N$ for which $\frac{\delta_{m-1}-\delta_{m}}{n_{m}}=:\ell_{m}\leq\tilde{\ell}_{m}$ holds define $f$ on $M_{m}$ by 
$$
f(x):=\left\{
\begin{array}{ll}
h(x_{k}\pm\delta_{m-1}\mp j\ell_{m})\  & \mbox{if }\, x=x_{k}\pm\delta_{m-1}\mp j\ell_{m}\\
& 0\leq j\leq n_{m},\; k=1,\dots,n\\
\mbox{affine on} & [x_{k}\pm\delta_{m-1}\mp j\ell_{m},x_{k}\pm\delta_{m-1}\mp (j+1)\ell_{m}]\\
& 0\leq j\leq n_{m}-1,\; k=1,\dots,n.
\end{array}
\right.
$$  
Finally set $f(x_{k}):=h(x_{k}),$ $k=1,\dots,n,$ such that by definition $f$ is continuous on $[0,1]$ and satisfies the $\lip$-condition for every point in $M:=\bigcup_{m=0}^{\infty}M_{m}=[0,1]\backslash\{x_{k}\}_{k=1}^{n}$ because it is locally polygonal around any such point.
 
We denote by $\overline{x}$ (resp.\ $\underline{x}$) the smallest (resp.\ biggest) element of $M$ having the form $x_{k}+\delta_{0}+j\ell_{0}^{k}$ or $x_{k}\pm\delta_{m-1}\mp j\ell_{m}$ which is bigger (resp.\ smaller) than or equal to $x$. For $x\neq \overline{x}$ an application of Lemma~\ref{KP} and our choice of $\ell^{k}_{0}$ and $\ell_{m}$ provide $L_{x\overline{x}}(f)\leq L_{\underline{x}\overline{x}}(h)\leq\frac{\eps}{4}$. A similar argument applies for $x\neq\underline{x}$ so that in the case $x,y\in M, x<y$, $x<\overline{x}<\underline{y}<y$ we obtain
\begin{align*}
L_{xy}(h-f) & \leq L_{x\overline{x}}(h-f)+L_{\overline{x}\underline{y}}(h-f)+L_{\underline{y}y}(h-f)\\
& \leq L_{x\overline{x}}(h)+L_{x\overline{x}}(f)+0+L_{\underline{y}y}(h)+L_{\underline{y}y}(f)\leq 4\cdot\frac{\eps}{4}=\eps.
\end{align*}
The other cases are treated similarly. Altogether we obtain $L_{\alpha}(h-f)\leq\eps$ and $f\in H_{\alpha}^{p}$ and conclude that $H_{\alpha}^{p}$ is dense in $H_{\alpha}^{\omega}$.
\end{proof}
 
We remark that by a combination of the above proofs one could obtain directly that $H_{\alpha}^{0}$ is an $M$-ideal in $H_{\alpha}^{\omega}$. Finally we answer the question whether $H_{\alpha}^{0}$ is even an $M$-summand and thus a trivial $M$-ideal in $H_{\alpha}^{\omega}$ negatively. Therefore we use (the easy part of) the following characterization of $M$-summands (see \cite[II.3.4]{HWW}) where we understand by $B(x,r)$ the closed ball around $x$ with radius $r$. 

\begin{prop}
\label{Msum}
A closed subspace $U$ of a Banach space $X$ is an $M$-summand in $X$ if and only if all families $\{B(x_{i},r_{i})\}_{i\in I}$ of closed balls with 
\begin{equation}
B(x_{i},r_{i})\cap U\neq\varnothing\quad\forall i\in I
\label{1infkp}
\end{equation}
and
\begin{equation}
\bigcap_{i\in I} B(x_{i},r_{i})\neq\varnothing
\label{2infkp}
\end{equation}
satisfy
\begin{equation*}
\bigcap_{i\in I} B(x_{i},r_{i})\cap U\neq\varnothing.
\end{equation*}
\end{prop}

Note that for $I=\{1,2,3\}$ the given condition requires the 3-ball property (in Proposition~\ref{3KP}) to be valid with $\eps=0$. 

\begin{prop}
\label{no Ms}
$H_{\alpha}^{0}$ is not an $M$-summand in $H_{\alpha}^{\omega}$.
\end{prop}

\begin{proof}
We show that the intersection condition given in Proposition~\ref{Msum} fails for $I=\{1,2\}$. Choose $h:x\mapsto x^{\alpha}$ in $B_{H_{\alpha}^{\omega}}$ and the functions $f_{1}:x\mapsto x$ and $f_{2}=-f_{1}:x\mapsto -x$ in $B_{H_{\alpha}^{0}}$. Then obviously
$$
f_{i}\in B(h+f_{i},1)\cap H_{\alpha}^{0}\quad\forall i\in\{1,2\}
$$
and
$$
h\in\bigcap_{i\in\{1,2\}} B(h+f_{i},1)
$$
holds thus conditions~\reff{1infkp} and~\reff{2infkp} are fulfilled. On the other hand there is no $g\in H_{\alpha}^{0}$ satisfying
\begin{equation}
L_{\alpha}(h+f_{i}-g)\leq 1\quad\forall i\in\{1,2\},
\label{no eps}
\end{equation}
and consequently
$$
\bigcap_{i\in\{1,2\}} B(h+f_{i},1)\cap H_{\alpha}^{0}=\varnothing.
$$ 
Assume  there is a $g\in H_{\alpha}^{0}$ with property~\reff{no eps},
without loss of generality real-valued. Then we get
$L_{01}(h+f_{1}-g)=|2-g(1)|\leq 1$ and $L_{01}(h+f_{2}-g)=|g(1)|\leq
1$, thus $g(1)=1$. Furthermore there is a point $x_{0}>0$ such that
$f_{1}(x_{0})+g(x_{0})<h(x_{0})$, since otherwise the inequality 
$$
\frac{|f_{1}(x)+g(x)-(f_{1}(0)+g(0))|}{|x-0|^{\alpha}}\geq\frac{|h(x)|}{x^{\alpha}}=1\quad\forall x\in (0,1]
$$
contradicts $f_{1}+g\in H_{\alpha}^{0}$. Now due to $h(1)-f_{1}(1)-g(1)=-1$ the mean value theorem provides a point $\tilde{x}\in (x_{0},1)$ with $h(\tilde{x})-f_{1}(\tilde{x})-g(\tilde{x})=0$ and we conclude
$$
L_{\tilde{x}1}(h+f_{2}-g)=\frac{|(h-f_{1}-g)(\tilde{x})-(h-f_{1}-g)(1))|}{|1-\tilde{x}|^{\alpha}}=\frac{1}{|1-\tilde{x}|^{\alpha}}>1,
$$ 
which contradicts~\reff{no eps}.
\end{proof}

Now collecting the facts stated in Propositions~\ref{Ha0 3B   Hap},
\ref{Hap dense Haw} 
and~\ref{no Ms} we obtain our main result.

\begin{thm}
\label{Ha0 pMid Haw}
$H_{\alpha}^{0}$ is a proper $M$-ideal in $H_{\alpha}^{\omega}$.
\end{thm}

\begin{rem}
\label{Mid-rem}
It is easy to see that the propositions leading to Theorem~\ref{Ha0 pMid Haw} are equally true for $\lip([0,1]^{\alpha})$ with the norm $\|f\|_{L}$ if one replaces $H^{p}_{\alpha}$ and $H^{\omega}_{\alpha}$ by $\Lip^{p}([0,1]^{\alpha})$ and $\Lip^{\omega}([0,1]^{\alpha})$ which are defined in the same manner.

Proposition~\ref{no Ms} states that there is no $M$-projection from
$H^{\omega}_{\alpha}$ onto $H^{0}_{\alpha}$. In fact $H^{0}_{\alpha}$
is not even complemented in $H^{\omega}_{\alpha}$. This follows from a
proof of Johnson in \cite[pp.~179--183]{John74} 
where, for metric spaces $K$ with a cluster point, he gives an isomorphic embedding $T:\ell^{\infty}\to \Lip(K^{\alpha})$ such that $T(c_{0})$ is complemented in $\lip(K^{\alpha})$. Now the construction of $T$ relies on the existence of one critical point only, thus we have $T(\ell^{\infty})\subset \Lip^{p}(K^{\alpha})$ which implies that in these general cases $\lip(K^{\alpha})$ is not complemented in $\Lip^{p}(K^{\alpha})$.  

More generally one can ask the following question. Taking into account
 the shifts of the function $x\mapsto x^{\alpha}$ one observes that
 $H^{\omega}_{\alpha}$ is nonseparable and moreover lies between
 $H^{0}_{\alpha}$, which is isomorphic to $c_{0}$, and $H_{\alpha}$,
 which is isomorphic to $\ell^{\infty}$. Can one derive just with this
 knowledge that $H^{0}_{\alpha}$ is not complemented in
 $H^{\omega}_{\alpha}$? The answer to this question is no. David Yost 
showed us a counterexample which we sketch here with his permission. The space $\ell^{\infty}$ contains a subspace $V$ isometric to $\ell^1(2^{\N})$ (see \cite[p.~155]{HHZ}); take $X$ to be the closed linear span of $c_{0}$ and $V$. Since the intersection $V\cap c_{0}$ is finite-dimensional (for instance, because $V$ has the Schur property), $c_{0}$ is complemented in the nonseparable space $X$.  
\end{rem}


\section{Generation of $H^{p}_{\alpha}$ and $H^{\omega}_{\alpha}$ using Ciesielski's basis}
\label{3}

Using the isomorphism $T: \ell^{\infty}\to H_{\alpha}$ which was found by Ciesielski in \cite{Cies60} we now want to give a description of $H^{p}_{\alpha}$ and $H^{\omega}_{\alpha}$ in terms of bounded sequences. It was already mentioned in the introduction that $T$ is defined by 
$$
T((a_{n}))=\sum_{n=1}^{\infty}a_{n}\varphi_{n},
$$
where $(\varphi_{n})$ is the Schauder basis for the space of all
continuous functions on $[0,1]$ vanishing at 0, normalized in
$H_{\alpha}$. For the construction of $(\varphi_{n})$ it is helpful
to decompose a natural number $n\geq 2$ as $n=2^{m}+k$ with
$m=0,1,\dots$ and $k=1,\dots,2^{m}$, where $m$ and $k$ are uniquely
determined by $n$. Except for $\varphi_{1}:x\mapsto x$ the graph of
any element in $(\varphi_{n})$ is an equilateral triangle on the
support $[\frac{k-1}{2^m},\frac{k}{2^m}]$ of $\varphi_{n}$. We
abbreviate $x_{n}^{\ell}:= \frac{k-1}{2^m}$, $x_{n}^{r}:=
\frac{k}{2^m}$ and $x_{n}^{c}:=(x_{n}^{r}+x_{n}^{\ell})/2$. 
It turns out that the inverse $T^{-1}$ is given by 
\begin{equation}
\label{ans}
T^{-1}(f)=\left(\frac{1}{2}\left(\frac{f(x_{n}^{c})-f(x_{n}^{\ell})}{(x_{n}^{c}-x_{n}^{\ell})^{\alpha}}-\frac{f(x_{n}^{r})-f(x_{n}^{c})}{(x_{n}^{r}-x_{n}^{c})^{\alpha}}\right)\right)_{n\in\N}
\end{equation}
so that local H\"older slopes of an $f\in H_{\alpha}$ are encoded in
the sequence $T^{-1}(f)$. Via $T$ we can now translate the definitions
of $H^{p}_{\alpha}$ and $H^{\omega}_{\alpha}$ into 
the sequence setting.

\begin{definition}
\label{copw}
We define $c_{p}$ to be the subspace of $\ell^{\infty}$ containing all sequences $(a_{n})$ with the following property:

There are finitely many points $x_{1},\dots,x_{N}\in [0,1]$ such that
any $x\notin\{x_{k}\}_{k=1}^{N}$ has a neighbourhood $U_x$ for which
the subsequence $(a_{n(\ell)})_{\ell\in\N}$ of $(a_{n})$ with
$\{n(\ell)\}_{\ell\in\N}=\{n=2^{m}+k:\frac{k-1}{2^m},\frac{k}{2^m}\in
U_{x}\}$ is in $c_{0}$.

We define $c_{\omega}$ to be the subspace of $\ell^{\infty}$ containing all sequences $(a_{n})$ with the following property:

For any $\eps>0$ there are finitely many points $x_{1},\dots,x_{N}\in
[0,1]$ such that any $x\notin\{x_{k}\}_{k=1}^{N}$ has a neighbourhood
$U_x$ for which the subsequence $(a_{n(\ell)})_{\ell\in\N}$ of
$(a_{n})$ with
$\{n(\ell)\}_{\ell\in\N}=\{n=2^{m}+k:\frac{k-1}{2^m},\frac{k}{2^m}\in
U_{x}\}$ satisfies $\|(a_{n(\ell)}) \|_{\infty}\leq\eps$. 
\end{definition}

It follows immediately from the definition of $H^{p}_{\alpha}$ and
$H^{\omega}_{\alpha}$ together with \reff{ans} that we have
$H^{p}_{\alpha}\subset T(c_{p})$ and $H^{\omega}_{\alpha}\subset
T(c_{\omega})$. Conversely, for $f\in T(c_{p})$ (or $f\in
T(c_{\omega})$) one can consider
$\tilde{f}=\sum_{\ell=1}^{\infty}a_{n(\ell)}\varphi_{n(\ell)}$ for any
given $x$ and $U_x$. Then we obtain $\tilde{f}\in H_{\alpha}^{0}$ (or
$L_{\alpha}(\tilde{f})\leq \|T\|\eps$) and $f=\tilde{f}+p$ locally
around $x$ with a polygon $p$. This observation provides the converse 
inclusions.

\begin{prop}
$H^{p}_{\alpha}$ and $H^{\omega}_{\alpha}$ are the images of $c_{p}$
 and $c_{\omega}$ respectively
 under Ciesielski's isomorphism $T$.
\end{prop}

We now turn to a reformulation of the properties given in
Definition~\ref{copw} which may be more instructive than the latter
because it focuses on the situation around the (critical) points
$x_{1},\dots,x_{N}\in [0,1]$ rather than on their 
complement and at the same time emphasizes a combinatorial point of view. 

\begin{prop}
$c_{p}$ is the space of all bounded sequences $(a_{n})$ with the following property:

There is a decomposition of the index set $\N=N_{0}\cup N_{1}$ into two nonintersecting subsets $N_{0}$ and $N_{1}$, such that $(a_{n})_{n\in N_{0}}\in c_{0}$ and $\{\frac{k}{2^{m}}:n=2^{m}+k\in N_{1}\}$ has finitely many cluster points.

$c_{\omega}$ is the space of all bounded sequences $(a_{n})$ with the following property:

For any $\eps>0$ there is a decomposition of the index set $\N=N_{0}\cup N_{1}$ into two nonintersecting subsets $N_{0}$ and $N_{1}$, such that $\|(a_{n})_{n\in N_{0}}\|_{\infty}\leq\eps$ and $\{\frac{k}{2^{m}}:n=2^{m}+k\in N_{1}\}$ has finitely many cluster points.
\end{prop}

\begin{proof}
That the sets with the properties given in the proposition are contained in $c_{p}$ or in $c_{\omega}$, respectively, is quite clear. 

So let $(a_{n})\in c_{p}$ and $x_{1},\dots,x_{N}\in [0,1]$ be given by Definition~\ref{copw}. Consider the set $C$ of all points $x\in [0,1]$ for which there exists a subsequence $(a_{n(\ell)})_{\ell\in\N}\notin c_{0}$ of $(a_{n})$ ($n(\ell)=2^{m(\ell)}+k(\ell)$) with $\frac{k(\ell)}{2^{m(\ell)}}\to x$ for $\ell\to\infty$. Choose these subsequences $(a_{n(\ell)})$ maximal in the sense that there is no non-$c_{0}$-subsequence of $(a_{n})$ disjoint to an $(a_{n(\ell)})$ and having the same limit as $(a_{n(\ell)})$ (this is an application of Zorn's Lemma in the countable setting). Now collect all indices $n(\ell)$ occurring in these subsequences in $N_{1}$. 

First it is obvious that $C$ is contained in $\{x_{k}\}_{k=1}^{N}$. We still have to show $(a_{n})_{n\in \N\backslash N_{0}}\in c_{0}$. Assume this is not the case; then there exists a subsequence $(a_{n(\ell)})_{\ell\in\N}$ of $(a_{n})_{n\in \N\backslash N_{0}}$ not belonging to $c_{0}$. On the other hand there is a cluster point $x\in [0,1]$ of the corresponding set $\{\frac{k(\ell)}{2^{m(\ell)}}\}_{\ell\in\N}$. We obtain $x\in C$ by definition, but this contradicts the maximality of $N_{1}$.   

The assertion for $c_{\omega}$ is proved in the same manner.
\end{proof}

\begin{rem}
\label{compare}
Observe that in contrast to $H^{p}_{\alpha}$ and $H^{\omega}_{\alpha}$
it is very easy to identify $c_{\omega}$ as the closure of $c_{p}$. On
the other hand if one translates the natural way to do this via $T$
into the function space $H^{\omega}_{\alpha}$ this leads to similar
``infinite polygonal'' approximations as considered in the proof of
Proposition~\ref{Hap dense Haw}. In a similar vein one can investigate
the method of the ``inserted constants'' as used for the proof of
Proposition~\ref{Ha0 3B Hap} when one formulates it in $c_{p}$. For a
given $h\in H_{\alpha}^{p}$ the approximation $g$ was defined to be
constant on intervals $I_{k}:=[x_{k}-\delta,x_{k}+\delta]$ around the
critical points $x_{k}$ of $h$ and elsewhere differs from $h$ only by
suitable constants. From~\reff{ans} we can see that the entries
$b_{n}$ in $T^{-1}(g)$ vanish if $\supp\varphi_{n}\subset I_{k}$ and
are equal to the entries $a_{n}$ in $T^{-1}(h)$ if
$\supp\varphi_{n}\subset [x_{k}+\delta,x_{k+1}-\delta]$ for a certain
$k$. This seems to be quite natural in $c_{p}$, too. However the
entries $a_{n}$ of $T^{-1}(h)$ with small indices have to be changed
according to the behaviour of $h$ around the critical points lying in
$\supp\varphi_{n}$ in order to get those of $T^{-1}(g)$. Of course
there would be no reason to do something like this in $c_{p}$, but for
the proof of the 3-ball property in $H^{p}_{\alpha}$ additive changes by polygons are essential.  
\end{rem}


\section{The almond function}
\label{4}

\begin{figure}[h]
\begin{minipage}[b]{0.49\linewidth}
\psfrag{1}{$1$}
\centering\epsfig{file=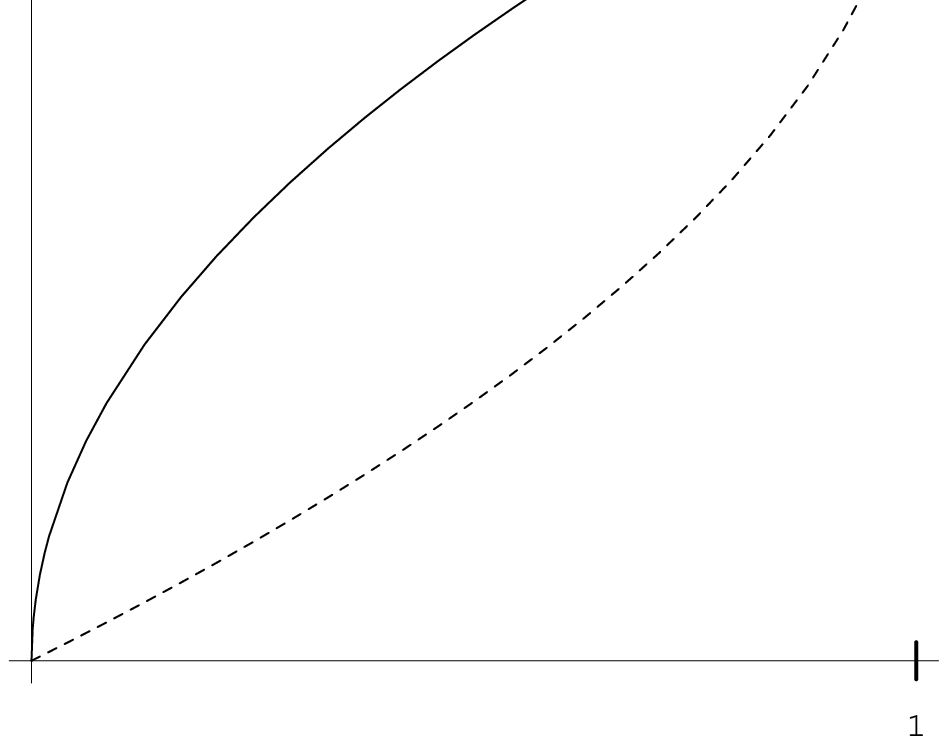,width=6cm}
\caption{$h_{0}$ and $\tilde{h}_{0}$}
\label{exman1}
\end{minipage}\hfill
\begin{minipage}[b]{0.49\linewidth}
\psfrag{r}{$r$}
\psfrag{s}{$s$}
\psfrag{1}{$1$}
\centering\epsfig{file=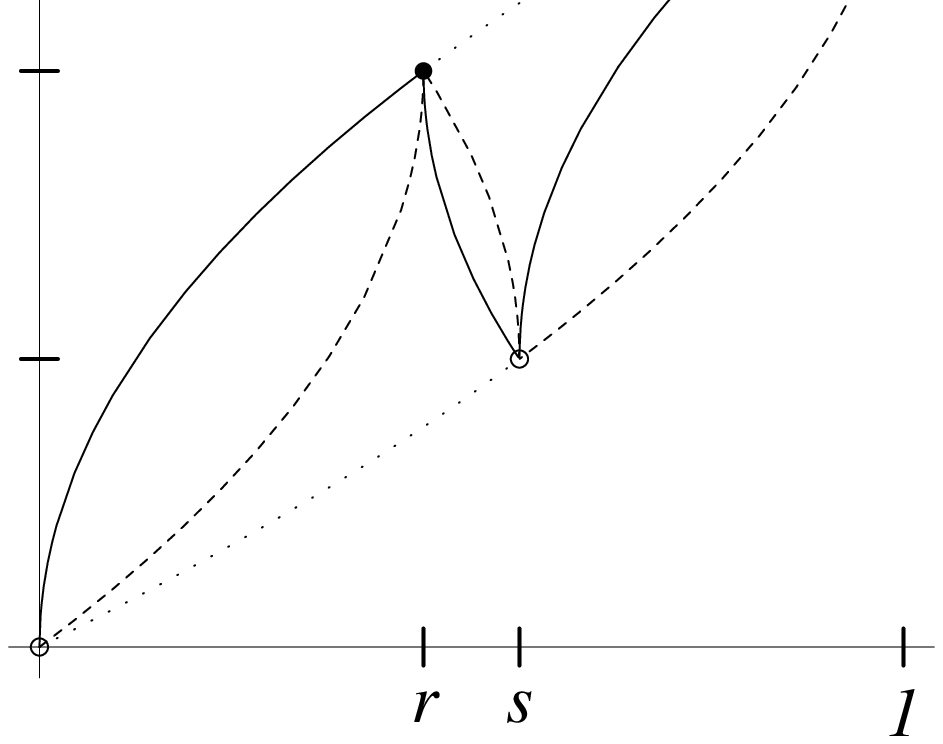,width=6cm}
\caption{First step}
\label{exman2}
\end{minipage}
\end{figure}

We now want to present an example of a function $h\in H_{\alpha}\backslash H^{\omega}_{\alpha}$ which is severely ``big Lipschitz'' in the sense that it has critical points everywhere. Consider the functions $h_{0}: x\mapsto x^{\alpha}$ and $\tilde{h}_{0}: x\mapsto 1-(1-x)^{\alpha}$. Since the area between the graphs of these two functions is shaped like an almond (see figure~\ref{exman1}) we simply call it an \emph{almond} for the sake of this construction. We find a point $r\in (0,1)$ such that the point $s\in (r,1)$ in which the graphs of $\tilde{h}_{0}$ and the function $x\mapsto r^{\alpha}-(x-r)^{\alpha}$ intersect satisfies $s=1-r$. Geometrically this results in cutting the first almond into three new ones the left and the right one of which are equal in size and look like the original almond. The boundary of the new almonds is given by the graphs of the two functions 

\begin{figure}[h]
\begin{minipage}[b]{0.49\linewidth}
\psfrag{r}{$r$}
\psfrag{s}{$s$}
\psfrag{1}{$1$}
\psfrag{p}{$r'$}
\psfrag{q}{$s'$}
\psfrag{u}{$t^{\alpha}$}
\psfrag{t}{$t$}
\centering\epsfig{file=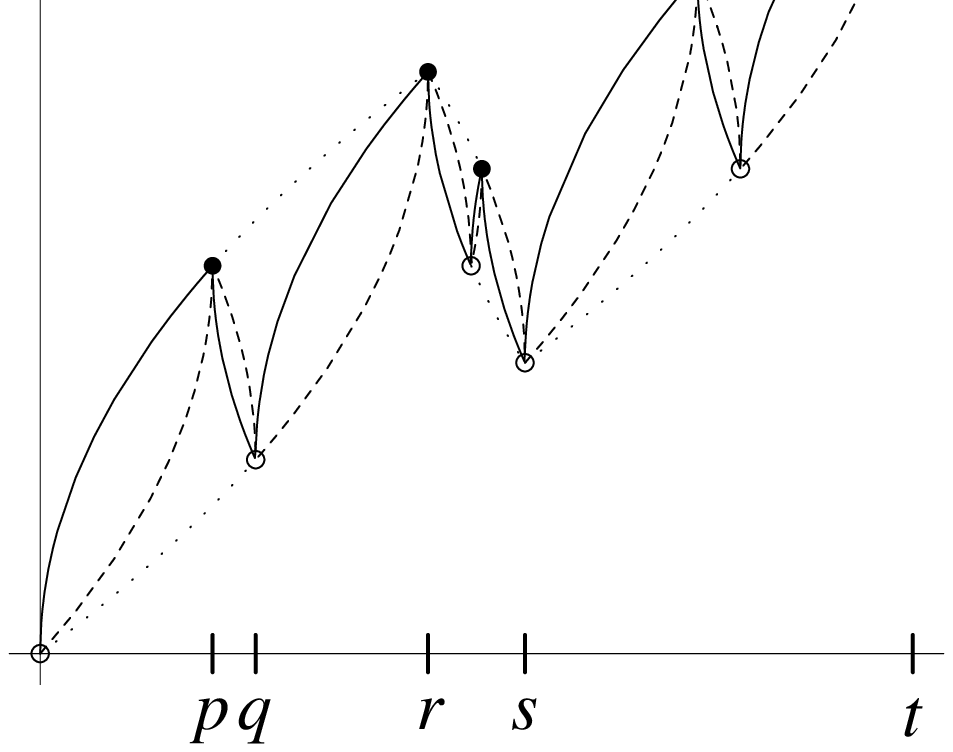,width=6cm}
\caption{Further construction}
\label{exman3}
\end{minipage}\hfill
\begin{minipage}[b]{0.49\linewidth}
\psfrag{1}{$1$}
\centering\epsfig{file=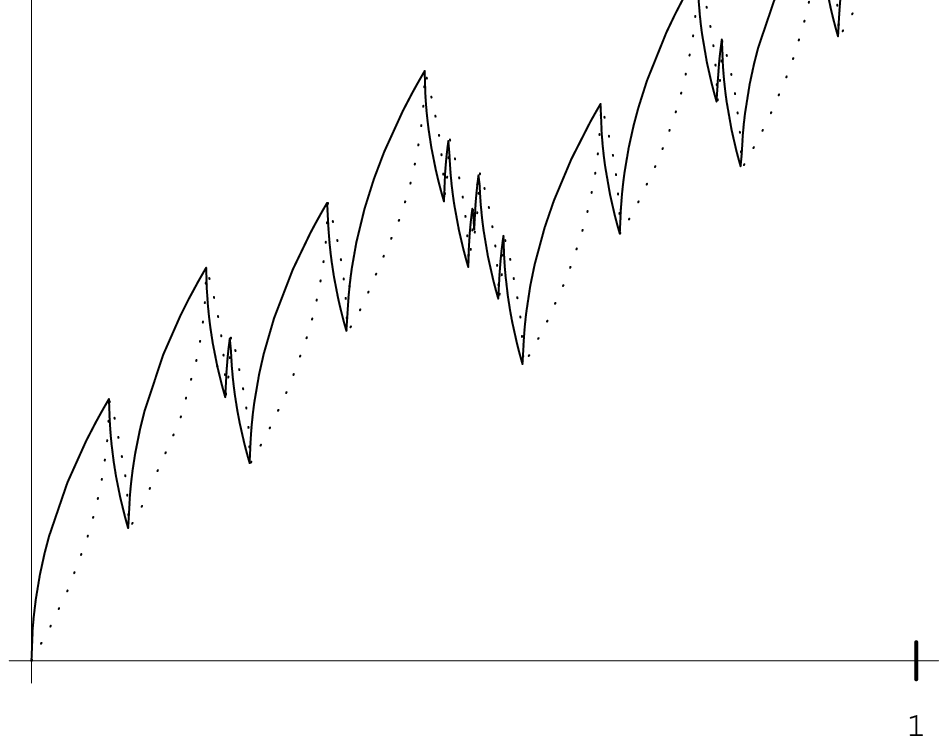,width=6cm}
\caption{$h_{3}$ and $\tilde{h}_{3}$}
\label{exman4}
\end{minipage}
\end{figure}

$$
h_{1}(x)=\left\{
\begin{array}{cl}
x^{\alpha} & \enspace 0\leq x\leq r\\
r^{\alpha}-(x-r)^{\alpha} & \enspace r\leq x\leq s\\
1-(1-s)^{\alpha}+(x-s)^{\alpha} & \enspace s\leq x\leq 1
\end{array}
\right.
$$
and
$$
\tilde{h}_{1}(x)=\left\{
\begin{array}{cl}
r^{\alpha}-(r-x)^{\alpha} & \enspace 0\leq x\leq r\\
1-(1-s)^{\alpha}+(s-x)^{\alpha} & \enspace r\leq x\leq s\\
1-(1-x)^{\alpha} & \enspace s\leq x\leq 1
\end{array}
\right.
$$
which can be seen in figure~\ref{exman2}. Now we can go on with this
procedure in the following way: If the orthogonal projection of an
almond onto the $x$-axis is of length $t$ we cut it as above into three
new almonds such that the orthogonal projections of the left and the
right one onto the $x$-axis both are of length $r$. 
This length is determined by the equation
$$
r^{\alpha}-(x-r)^{\alpha}=t^{\alpha}-(t-x)^{\alpha}
$$
where $x=t-r$ such that the ratio $y=r/t\in (0,1/2)$ which is given by 
\begin{equation}
\label{ka}
2y^{\alpha}-1=(1-2y)^{\alpha}
\end{equation}
is a monotonic increasing function $k$ of $\alpha\in (0,1)$. For example one obtains $k(1/2)=4/9$. Now it is easy to define suitable functions $h_{2}$ and $\tilde{h}_{2}$ whose graphs are the boundary of the nine almonds obtained from cutting the three almonds above. We illustrate suitable definitions in figure~\ref{exman3} (with $t=1$) and figure~\ref{exman4} for the second and the third step.

Going on with this construction one can get a sequence of functions $(h_{n})$ and $(\tilde{h}_{n})$ converging uniformly to the same limit.

\begin{definition}
\label{almond function}
The uniform limit $h$ of the sequence $(h_{n})$ is called \emph{almond function}.
\end{definition}

It is easy to see that like the functions $h_{n}$ and $\tilde{h}_{n}$
the almond function is an extreme point in the unit ball of
$H_{\alpha}$. Notice that the almond function is already uniquely
determined by the top points (see $\bullet$ in figure~\ref{exman2} and
\ref{exman3}) and the bottom points (see $\circ$ in
figure~\ref{exman2} and~\ref{exman3}) of the almonds appearing in the
construction. Once these points have occurred they remain untouched in
the following steps, i.e., any local minimum (maximum) of an $h_{k}$
is a local minimum (maximum) of any $h_{n}$ with $n\geq k$ and thus a
local minimum (maximum) of the almond function. But more can be said
about the properties of $h$ on the set $M$ of the abcisses of all such
points. By approaching maxima of $h$ by mimima and vice versa it
becomes clear from the construction that $h$ behaves quite badly 
on the dense subset $M$ of $[0,1]$.

\begin{prop}
\label{bad almond function}
The almond function $h$ is an extreme point of $B_{H_{\alpha}}$ and satisfies
$$
\limsup_{\substack{y\to x\\ x\neq y}} L_{xy}(h)=1\quad\forall x\in M.
$$
In particular $h$ only has critical points.
\end{prop}

Now approaching minima (resp.\ maxima) of $h$ by minima (resp.\ maxima) even provides a positive lower bound for the slopes of $h$ around the points in $M$. In the computation the function $k:\alpha\mapsto y$ defined by the equation in~\reff{ka} is helpful. 

\begin{prop}
\label{very bad almond function}
The almond function $h$ satisfies
$$
\liminf_{\substack{y\to x\\ x\neq y}} L_{xy}(h)=\frac{1-k(\alpha)^{\alpha}}{(1-k(\alpha))^{\alpha}}\quad\forall x\in M.
$$
\end{prop}

\begin{proof}
It is enough to analyse the behaviour of $h$ at its minima around the point $0$. A look at figure~\ref{exman3} makes it clear that $h(s)=t^{\alpha}-(k(\alpha)t)^{\alpha}$ for $s=(1-k(\alpha))t$. Consequently the minima of $h$ lie above the graph of the function $f:x\mapsto (1-k(\alpha)^{\alpha})(1-k(\alpha))^{-\alpha}x^{\alpha}$ and it is easily checked that $f\geq h$ (locally) around $0$ with equality for the type of minima $(s,h(s))$ just considered. 
\end{proof}

The Propositions~\ref{bad almond function} and~\ref{very bad almond function} suggest that the almond function might be a sufficiently ``bad'' big Lipschitz function for our purposes. And in fact we can prove that the straightforward idea to approximate a big Lipschitz function $h$ by a polygon interpolating $h$ in its nodes in order to get the 3-ball property fulfilled via Lemma~\ref{3B} won't work in this example, however fine one chooses the step size. 

\begin{prop}
\label{quite bad almond function}
There is a constant $c>0$ depending only on $\alpha$ such that for any polygon $g$ interpolating the almond function $h$ in its nodes the estimate
$$
L_{xy}(h-g)\geq 1+c
$$
holds true for certain $x,y\in (0,\tilde{x})$ where $\tilde{x}$ is the smallest positive number with $g(\tilde{x})=h(\tilde{x})$. 
\end{prop}

\begin{proof}
Let $0=x_{0}<x_{1}<\dots<x_{n}=1$ be the partition of $[0,1]$ in which
a polygon $g$ has its nodes and let $\tilde{x}\leq x_{1}$ be the
smallest positive number with $g(\tilde{x})=h(\tilde{x})$. It is clear
from the construction of the almond function $h$ that there exists a
$t>0$ at which $h$ has a maximum such that with $r=k(\alpha)t$ we have
$\tilde{x}\in [r,t]$ and we are in the situation given in figure~\ref{exman3}
with $s=t-r$, $r'=k(\alpha)r$ and $s'=r-r'$. Then of course 
we obtain
$$
\frac{g(x)-g(y)}{x-y}\geq\frac{h(s)}{s}=\frac{1-k(\alpha)^{\alpha}}{1-k(\alpha)}t^{\alpha-1}\quad\forall
x,y\in [0,\tilde{x}], \ x\neq y.
$$
In particular this leads to
$$
\frac{g(s')-g(r')}{(s'-r')^{\alpha}}
\geq\frac{(1-k(\alpha)^{\alpha}) (k(\alpha)(1-2k(\alpha)))^{1-\alpha}}
{1-k(\alpha)}=:c>0
$$
and since $\frac{h(s')-h(r')}{(s'-r')^{\alpha}}=-1$ the choice $x=r'$ and $y=s'$ proves the assertion.
\end{proof}

The proof of this proposition makes it clear that the criterion for
the $3$-ball property given in Lemma~\ref{3B} is not satisfied in the case of the almond function $h$ for any approximating polygon $g$ intersecting $h$ between its first two nodes. The reason for this is that $g$ would be too steep on a too big interval $[0,\tilde{x}]$. Of course, in order to apply Lemma~\ref{3B}, $g$ could still be chosen that steep on a smaller interval since $g\in H_{\alpha}^{0}$. One is tempted to construct $g$ around 0 in such a way that it ``respects'' big ``dangerous'' intervals (such as $[r',s']$ in the proof) being more or less constant on them and still ``gains height'' on the complement of these intervals where one makes use of the $\eps>0$ that we have at our disposal for proving the $3$-ball property. In order to get such an approximation it seems necessary to look closely at the quantiative properties of the almond function on quite different scales. We don't know if such an approximation exists.

\begin{rem}
A look at Ciesielski's isomorphism $T: \ell^{\infty}\to H_{\alpha}$
discussed in 
Section~\ref{3} suggests that one might also obtain a quite extreme
element of $H_{\alpha}$ by a more analytic approach considering
$f=T((1,1,1, \dots))$. And indeed, by calculating a lower bound $b>0$
for $\|T\|$ Ciesielski himself proves that $f$ has only critical
points. More precisely he considers $f_{N}=\sum_{n=1}^{N}\varphi_{n}$
for $N=2^{m}+1$ and a certain positive null sequence $(x_{j})$ 
for which an equality 
$$
L_{0x_{j}}(f_{N})=b-\eps_{1}(j)-\eps_{2}(j,N)
$$
holds where $\lim_{j\to\infty}\eps_{1}(j)=0$ and
$\lim_{N\to\infty}\eps_{2}(j,N)=0$ for any $j\in\N$. Herefrom one
concludes $\lim_{j\to\infty}L_{0x_{j}}(f)=b$. Moreover the situation
in any point $y\in M':=\{\frac{k}{2^{m}}:n=2^{m}+k\in\N\}$ is in
principle the same as in 0 (up to a polygon which is in
$H_{\alpha}^{0}$) so that $\lim_{z_{j}\to
  0}L_{0z_{j}}(f)=\lim_{y+z_{j}\to y+}L_{y,y+z_{j}}(f)$ for any $y\in
M'\backslash\{1\}$ and any positive null sequence $(z_{j})$. So at
least in the sense of Proposition~\ref{bad almond function} the
function $f$ has similar properties as $h$. Still we find that our
geometric construction of the almond function (which for $\alpha=1/2$
was already considered in another 
context, see \cite[pp.~22]{Man}) might be more helpful for our purposes. 

Also note that as already discussed in Remark~\ref{compare} the
natural way to approximate $H_{\alpha}$-functions by
$H_{\alpha}^{0}$-functions in the sense of Ciesielski's isomorphism
ignores differences of the slope of two functions considered. We point
out that the idea to approximate an $h\in H_{\alpha}$ by the polygon
$g$ ``built up'' by the first $N$ entries of $T^{-1}(h)$ leads exactly
to the problems demonstrated in Proposition~\ref{quite bad almond
  function} since $g$ would then be a polygon interpolating $h$ in all
its nodes which are the centers of all intervals $\supp\varphi_{n}$
for $n\leq N$. So at this stage the question whether the almond
function is a counterexample to the assertion ``$H_{\alpha}^{0}$ is an
$M$-ideal in $H_{\alpha}$'' is still 
open.
\end{rem}


\end{document}